\documentclass[11pt]{article}

\usepackage{amsfonts,amsmath,amssymb,graphics}

\newcommand{\C}{\mathbb{C}}

\newcommand{\I}{\mathrm{i}}
\newcommand{\pr}{\mathbf{P}}
\newcommand{\ex}{\mathbf{E}}

\newcommand{\hmu}{\overline{\mu}}
\newcommand{\hsigma}{\overline{\sigma}}
\newcommand{\hkappa}{\overline{\kappa}}

\newcommand{\notthis}[1]{}

\newcommand{\Imag}{\mathrm{Im}\,}
\newcommand{\inv}{^{-1}}
\newcommand{\indic}[1]{\mathbf{1}_{#1}}

\newcommand{\invintr}{\int_{-\infty}^\infty}
\newcommand{\invintip}{\int_{-\I\infty}^{\I\infty\ast}}

\newcommand{\invintrp}{\int_{-\infty}^{\infty\ast}}

\newcommand{\half}{\frac{1}{2}}

\newcommand{\cdl}{\,|\,}

\newcommand{\shalf}{{\textstyle\frac{1}{2}}}

\newcommand{\pderiv}[2]{\frac{\partial{#1}}{\partial{#2}}}
\newcommand{\deriv}[2]{\frac{d{#1}}{d{#2}}}

\begin{document}

\title{\bf Time since maximum of Brownian motion and asymmetric L\'evy processes}

\author{R. J. Martin\footnote{Department of Mathematics, Imperial College London, London SW7 2AZ, U.K.}\mbox{ } and M. J. Kearney\footnote{Senate House, University of Surrey, Guildford GU2 7XH, U.K.}}
\maketitle

\begin{abstract}

Motivated by recent studies of record statistics in relation to strongly correlated time series, we consider explicitly the drawdown time of a L\'evy process, which is defined as the time since it last achieved its running maximum when observed over a fixed time period $[0,T]$. We show that the density function of this drawdown time, in the case of a completely asymmetric jump process, may be factored as a function of $t$ multiplied by a function of $T-t$. This extends a known result for the case of pure Brownian motion. We state the factors explicitly for the cases of exponential down-jumps with drift, and for the downward Inverse Gaussian L\'evy process with drift.

Published in J.\ Phys.\ A (Math.\ Theor.), 51:275001, 2018 (this version has minor updates)
%AMS 60G51, 60G17

\end{abstract}

\section{Introduction}

Interest in the statistics of records, closely linked to studies of extreme statistics, has rekindled over the last few years with significant progress made in terms of developing new theoretical insights. Early work was concerned with discrete processes where the variables are independent and identically distributed (i.i.d; for a review see \cite{Wergen13}). However, attention has recently shifted to consider records in the context of strongly correlated time series and in particular those whose \emph{increments} are independent, as in a Brownian motion. The exact results obtained in \cite{Majumdar08b} in relation to random walks with a continuous and symmetric jump distribution, concerning the number of records in a given interval as well as details of their duration, has led to many parallel studies which have greatly expanded our understanding. This includes the extension to random walk sequences with drift \cite{LeDoussal09,Majumdar12,Gouet15}, to continuous-time random walks \cite{Sabhapandit11}, and deeper examination of the special role played by the shortest, longest and last record in a given interval \cite{Godreche09,Godreche14}. An excellent overview of all this work and the plethora of areas of application in the physical and social sciences may be found in \cite{Godreche17}.

The question as to the duration of the last record is of particular interest in the context of financial time series, where more general studies of record processes \cite{Wergen11,Sabir14} have been augmented by consideration of this specific `time since maximum' or `drawdown time' \cite{Mijatovic12,Challet15,Challet17}. This is one of the motivations for the current work. However, we are also interested in how to generalise the results that have been obtained in relation to random walk or L\'evy processes. A great deal of current understanding (see e.g.~\cite[\S3]{Godreche17}) is based on a fundamental application of a celebrated theorem of Sparre-Andersen \cite{SparreA53,SparreA55}, the utilisation of which is relatively straightforward if the jump distribution is symmetric, but otherwise much more challenging. Another motivation,  therefore, is to consider how to handle situations where in the context of a generic L\'evy process (including drift and diffusion) the jump distribution is totally asymmetric, i.e.\ jumps occur in one direction only.

% geology: deposition or erosion
% meteorol
% insurance maths
% queues, road traffic

By way of an introduction and to set further work in context, let us consider the arithmetic Brownian motion
\[
dX_t = \mu\, dt + \sigma\, dW_t
\]
and define the running maximum to be
\[
M_t=\max_{0\le t'\le t}X_{t'}.
\]
Over a time period $[0,T]$ we wish to study the time at which, or since, the maximum occurred.
We define the drawdown process to be $(M_t-X_t)$, and consider
\[
\tau = \sup_{t \le T} \{ T-t \; : \; M_T = M_t \}
\]
which is the drawdown time.
For the Brownian motion with drift, the p.d.f.\ of $\tau$ is well known \cite{Borodin02,Buffet03}:
\footnote{As usual $\phi$, $\Phi$ denote the density and cumulative of the standard Normal distribution, while $\Phi_2$, when it occurs later, is the cumulative bivariate Normal distribution. }
\begin{equation}
f_\tau(t) = \frac{2}{\sqrt{t(T-t)}} \,
\mathcal{C}\!\left(\frac{\mu\sqrt{T-t}}{\sigma}\right)
\mathcal{C}\!\left(\frac{-\mu\sqrt{t}}{\sigma}\right)
\indic{0<t<T}, \qquad
\mathcal{C}(x) \equiv \phi(x)+x\Phi(x),
\label{eq:taudens}
\end{equation}
and reduces to the arcsine law when there is no drift.
Its cumulative, on the other hand, is not well known, and to the best of our knowledge the following is a new result:
\begin{eqnarray}
\pr(\tau<t) &=& \Phi\!\left(\frac{\mu\sqrt{T}}{\sigma}\right) - \frac{\mu\sqrt{T}}{\sigma} \phi\!\left(\frac{\mu\sqrt{T}}{\sigma}\right) 
\label{eq:taucumul}
 \\
&& \mbox{}+2 \left(1+\frac{\mu^2(T-t)}{\sigma^2}\right) \Phi_2\!\left(\frac{\mu\sqrt{T-t}}{\sigma},\frac{-\mu\sqrt{T}}{\sigma};-\sqrt{\frac{T-t}{T}}\right) \nonumber \\
&& \mbox{}-2 \left(1+\frac{\mu^2t}{\sigma^2}\right) \Phi_2\!\left(\frac{-\mu\sqrt{t}}{\sigma},\frac{\mu\sqrt{T}}{\sigma};-\sqrt{\frac{t}{T}}\right) \nonumber \\
&& \mbox{} + \frac{2\mu\sqrt{T-t}}{\sigma} \, \phi\!\left(\frac{\mu\sqrt{T-t}}{\sigma}\right) \Phi\!\left(\frac{-\mu\sqrt{t}}{\sigma}\right) + \frac{2\mu\sqrt{t}}{\sigma} \, \phi\!\left(\frac{\mu\sqrt{t}}{\sigma}\right)  \Phi\!\left(\frac{\mu\sqrt{T-t}}{\sigma}\right) ,
\nonumber
\end{eqnarray}
which is considerably less elegant than (\ref{eq:taudens}) but nevertheless of comparable utility. If one goes about deriving these results `long-hand', i.e.\ with no sophisticated techniques, the derivation is quite messy, and for that reason given in the Appendix, along with various useful integrals needed en route\footnote{As will become apparent, a derivation by Wiener-Hopf is much quicker.}. One important tool used therein, and having wider ramifications, is what we are calling the \emph{reciprocity law} for the bivariate Normal integral, linking $\Phi_2(x,-\rho x;\rho)$ and $\Phi_2(-x,\rho^*x;\rho^*)$
where $\rho^*=\sqrt{1-\rho^2}$. Incidentally a byproduct of this is the intriguing identity
\[
\Phi_2\! \left( x,\frac{x}{\sqrt{2}} ; \frac{1}{\sqrt{2}}  \right) = \frac{1}{2} \left( \Phi\!\left( \frac{x}{\sqrt{2}}\right)^2 + \Phi(x) \right).
\]
The appearance of the bivariate Normal integral in (\ref{eq:taucumul}) is expected given the orthant law that links this to the arcsine function encountered when $\mu=0$. Interesting as they are, however, these points are not the main focus of the paper, which is to investigate more general dynamics than a simple Brownian motion.

It is worth mentioning here that transformation by a strictly monotone function, so that $Y_t=\psi(X_t)$, has no effect on the distribution of $\tau$, so that departure from an arithmetic Brownian motion does not necessarily vitiate (\ref{eq:taudens},\ref{eq:taucumul}). An obvious example of this is the geometric Brownian motion.
Nonetheless the Brownian motion is restrictive because its trajectories are continuous: certainly in the context of financial modelling, and in particular the modelling of credit risks, this assumption is unjustifiable (see e.g.~\cite{Martin09a,Martin10b} and references therein).

A natural way of relaxing this is to use a L\'evy process, thereby retaining i.i.d.\ increments but dropping the assumption that they be Gaussian.
The tails are fatter and are usually either exponential, in which case all moments exist but only some exponential moments, or alternatively power-law, in which case only some moments, and no exponential moments, exist. A useful way of thinking about a L\'evy process is as a Brownian motion plus an array of mutually independent jump processes, each of which produces jumps of a given size; these have different hazard rates, and the higher the rate of occurrence of large jumps the more fat-tailed the process. This is enshrined in the L\'evy-Khinchin theorem (see e.g.\ \cite{Sato02}).
If we wish to make more concrete statements about the distribution of $X_t$, the independence of additive increments implies a convolution of the probability densities, and this is most easily represented by the characteristic function, $\ex[e^{\I u X_t}]$. This must be of the form $e^{L(u)t}$ for some function $L$ which we call the L\'evy generator. At one level we can simply view L\'evy processes by reference to their generators, with the L\'evy-Khinchin theorem giving necessary and sufficient conditions for validity of $L$.

As a specific example, let us consider the possibility of replacing the diffusion term with an exponential jump process, i.e.\ a Poisson process of rate $\lambda$ with jumps exponentially distributed of mean $\xi$, while retaining the drift term.
In the limit of high jump intensity and small mean jump size, with $\lambda\xi$ held fixed, we recover the Brownian motion. So this model is a three-parameter model that has the Brownian motion as one special case, but allows for the third moment to be nonzero. Another recipe is a L\'evy process such as the Inverse Gaussian \cite{Schoutens03} which is another way of departing from the Brownian case at the expense of adding only one extra parameter, or the Carr-Madan-Yor model \cite{Carr04,Schoutens03} in which two extra parameters are used.

The difficulty with generalising the Brownian motion is that it is considerably harder to obtain analytical results.
However, one thing about (\ref{eq:taudens}) is striking: it factorises as the product of two related pieces, one dependent on $t$ and the other on $T-t$. The derivation in the Appendix does not make this at all obvious, and in fact it relies on a result concerning the bivariate Normal integral, which we call the reciprocity law, allowing the results to be tidied up. As we will presently show, a different method of derivation using the so-called fluctuation identities allows an analogue of (\ref{eq:taudens}) to be obtained for L\'evy processes---at least in the special case of jumps in one direction only, which is known as the \emph{asymmetric} L\'evy case. In principle this case can be analysed starting from \cite{Pitman81}, but this is far from obvious for a researcher lacking the necessary background in probability: the derivation given here is more straightforward, requiring essentially no more than the manipulation of Fourier and Laplace  integrals. Besides, integral representations are useful in the computation of the results.  
Another facet of (\ref{eq:taudens}) is worth commenting on: the form of $\mathcal{C}(x)$ suggests a correspondence\footnote{Cf.~the call and put payoffs in a Normal option-pricing model, or for a more general treatment 
\cite[\S2.5]{Martin11b}. }
with the expectations $\ex[\max(X_{T-t},0)]$ and $\ex[\max(-X_t,0)]$. A natural question to ask is whether this correspondence also carries over to the asymmetric L\'evy case.
As it happens, the answer to that is ``yes and no'': one of the factors does correspond in the suggested way but the other does not. These results are encapsulated by (\ref{eq:wh6}) \emph{et seq}.

The derivation of these results requires the fluctuation or Wiener-Hopf identities, which we discuss in the next section.
The exponential-jumps model with drift but no diffusion turns out to be reasonably tractable, mainly because its L\'evy generator is a simple algebraic function and the ensuing Fourier integrals invoke functions `no worse than' the confluent hypergeometric function. We can then obtain the density of $\tau$, via (\ref{eq:wh6}) and the two formulae (\ref{eq:edj1},\ref{eq:edj2}). Much the same remarks apply to our other example of choice, the Inverse Gaussian process plus drift.

\section{Analysis}

The Appendix gives a step-by-step derivation of the results (\ref{eq:taudens},\ref{eq:taucumul}) for the Brownian motion, and relies on the following argument: for the maximum to have occurred on or after time $t$, it is necessary and sufficient that $M_T>M_t$. In other words $\pr(M_T>M_t)=\pr(\tau<T-t)$.
To find this, we take the joint density of $(X_t,M_t)$ and the conditional distribution of the running maximum given a particular starting-point. We then integrate out.
Here we follow the same procedure but carry out most of the work in Fourier space, and using the fluctuation identities or Wiener-Hopf factorisation.

Let the process be described by its L\'evy generator $L(u)$, so that $\ex[e^{\I u X_t}] = e^{L(u)t}$.
For the arithmetic Brownian motion with drift, $L(u) = \I \mu u-\half \sigma^2 u^2 $.
For other processes one can express $L$ either directly or through its L\'evy-Khinchin representation\footnote{The version we state is simplified, and applies to processes of finite variation.} as a Brownian motion plus an array of jumps of all sizes, weighted by the L\'evy measure $\nu(x)\ge0$:
\begin{equation}
L(u) = \I \mu u - \shalf \sigma^2 u^2 + \int_{-\infty}^\infty (e^{\I u x} - 1) \nu(x)\, dx.
\label{eq:lkh}
\end{equation}
We recall (see e.g.~\cite{Sato02} for details) the fluctuation identity
\begin{equation}
\int_0^\infty s e^{-st} \ex[e^{\I vX_t+\I w M_t}] \, dt = \varphi^+_s(v+w) \varphi^-_s(v)
\label{eq:fluc1}
\end{equation}
where for each $s$ the functions $\varphi^+_s(\cdot)$, $\varphi^-_s(\cdot)$ are analytic and bounded and free from zeros in the upper and lower half-planes respectively, with
\begin{equation}
\varphi^+_s(u)\varphi^-_s(u) = \frac{s}{s-L(u)}, \qquad \varphi^+_s(0)=\varphi^-_s(0)= 1
\end{equation}
commonly known as the Wiener-Hopf factorisation.

To find $\pr(M_T>M_t)$, we first condition on the pair $(X_t,M_t)$.
We have as a direct consequence of (\ref{eq:fluc1})
\begin{equation}
f_{X_t,M_t}(x,y) = \frac{-1}{(2\pi\I)^3} \invintip \invintr\invintr \varphi^+_{s_1}(u+v) \varphi^-_{s_1}(u)  e^{-\I u x} e^{-\I v y} e^{s_1t}  \, du \, dv \, \frac{ds_1}{s_1}
\label{eq:wh1}
\end{equation}
(* denotes that the path of integration is deformed so as to pass around the origin anticlockwise)
and
\begin{equation}
\pr(M_t>z) = \frac{1}{(2\pi\I)^2} \invintrp \invintip \varphi^+_{s_1}(w) e^{s_1 t} e^{-\I w z} \, \frac{ds_1}{s_1} \, \frac{dw}{w}
\label{eq:wh2}
\end{equation}
so that
\begin{equation}
\pr \!\left( \max_{t'\in[t,T]} X_{t'} >y \,\big|\, X_t=x \right) = \frac{1}{(2\pi\I)^2} \invintrp \invintip \varphi^+_{s_2}(w) e^{s_2 (T-t)} e^{-\I w (y-x)} \, \frac{ds_2}{s_2} \, \frac{dw}{w}.
\label{eq:wh3}
\end{equation}
Now take $(\ref{eq:wh1})\times(\ref{eq:wh3})$ and integrate over $y\in[0,\infty)$ and $x\in(-\infty,y]$ which will give $\pr(M_T>M_t)$. Two pairs of integrals---$\int\int\,du\,dx$ and $\int\int\,dy\,dv$---collapse immediately (Fourier inversion), and on interchanging $t\leftrightarrow T-t$ we have
\begin{equation}
\pr(\tau<t) = \frac{1}{(2\pi\I)^3} \invintrp\invintip\invintip \varphi^-_{s_1}(w) \varphi^+_{s_2}(w) e^{s_1(T-t)+s_2t} \, \frac{ds_1}{s_1} \, \frac{ds_2}{s_2} \,  \frac{dw}{w} + \indic{t>T}.
\label{eq:wh4}
\end{equation}
On differentiating w.r.t.~$t$ we obtain
\begin{equation}
f_\tau(t)  = \frac{1}{(2\pi\I)^3} \invintrp\invintip\invintip \varphi^-_{s_1}(w) \varphi^+_{s_2}(w) e^{s_1(T-t)+s_2t} \, (s_2-s_1) \frac{ds_1}{s_1} \, \frac{ds_2}{s_2} \,  \frac{dw}{w}
+ \delta(T-t).
\label{eq:wh5}
\end{equation}
The next step has to be to attempt the $w$-integral. Note that in (\ref{eq:wh5}) the singularity at $w=0$ generates no contribution, provided $0<t<T$, as the integral is $\delta(t)\indic{t<T}-\delta(T-t)\indic{t>0}$. To make further progress we need to understand where the more important singularities lie in the $w$-plane.

From this point onwards our development is similar to that in \cite{Madan08}, who examine the related problem of a spectrally-negative (no up-jumps) process hitting a barrier; this and related papers stem from work done by Lipton \cite{Lipton02c}. The point about all this work is that one of the Wiener-Hopf factors is a simple pole. This can be seen by applying the Argument Principle to the function $\frac{s}{s-L(u)}$ around the real axis closed in the lower half-plane, and using (\ref{eq:lkh}) to figure out the behaviour of $L$. The position of this pole is easily found, and then the other factor can be inferred immediately. Unlike the analysis in \cite{Madan08}, however, we give explicit solutions for the inversion integrals in the case we consider, whereas they devise a `black-box' method that will analyse any process using numerical Laplace inversion. As a result we can deduce concrete results in these cases, admittedly invoking special functions.

% M and S provide a sort of black-box.

\subsection{Processes with no up-jumps}

When the process has no up-jumps the factor $\varphi^+$ has only one singularity in the lower half-plane, which is a simple pole, so that the factorisation is
\[
\varphi^+_s(u) = \frac{u^+_s}{u^+_s-u}, \qquad \varphi^-_s(u) = \frac{s}{s-L(u)} \biggr/ \frac{u^+_s}{u^+_s-u}, \qquad \Imag u^+_s<0.
\]
If we pull the $w$-contour down through the pole $w=u^+_s$ and out to $\infty$ (whereupon the integral vanishes), and evaluate the residue at $u^+_s$, we find the neat result
\[
f_\tau(t)  = \frac{1}{(2\pi\I)^2} \invintip\invintip \frac{u^+_{s_2}}{u^+_{s_1}} e^{s_1(T-t)+s_2t}  \, ds_1\, \frac{ds_2}{s_2} 
\]
(note $L(u^+_s)=s$, which explains some of the cancellation).
Accordingly
\begin{equation}
f_\tau(t)  = 
\biggr( \frac{1}{2\pi\I} \invintip \frac{1}{\I u^+_{s_1}} e^{s_1(T-t)} \, ds_1 \biggr) \biggr( \frac{1}{2\pi\I} \invintip \frac{\I u^+_{s_2}}{s_2} e^{s_2t}  \, ds_2 \biggr),
\label{eq:wh6}
\end{equation}
factorising as we had hoped. 
Changing variable from $s$ to $u$ where $L(u)=s$, we have
\begin{equation}
f_\tau(t)  = 
\biggr( \frac{-1}{2\pi} \invintrp \frac{L'(u_1)}{u_1} e^{L(u_1)(T-t)} \, du_1 \biggr)
\biggr( \frac{1}{2\pi} \int_\gamma \frac{u_2L'(u_2)}{L(u_2)} e^{L(u_2)t}  \, du_2 \biggr)
\label{eq:wh7}
\end{equation}
where the contour $\gamma$ runs from $-\infty$ to $\infty$, but care is needed over singularities (there is none at $u=0$, but there may be others).
Now if $C(u)=\ex[e^{iuX}]$ denotes the characteristic function of a random variable $X$, then
\[
\ex [ \max(X,0) ] = -\frac{1}{2\pi} \invintrp C(u) \, \frac{du}{u^2}
\]
so that the first factor can be integrated by parts to give $\ex[\max(X_{T-t},0)]/(T-t)$, the form of which was anticipated earlier.
But the second factor is not the same, and its interpretation is seen via (\ref{eq:wh2})---pulling the $w$-contour down through the pole, as before---to be the density of $M_t$ evaluated at the starting-point $X_0$ (=0), or 
\[
f_{M_t}(0) = \lim_{y\searrow0} \frac{\pr(M_t<y)}{y} .
\]
%\footnote{For an $\alpha$-stable process it is simply $\alpha f_X(0)$. However, that is unhelpful here.}. 
Thus
\begin{equation}
\mbox{ [No up-jumps] }\qquad
f_\tau(t)  = 
\frac{1}{T-t} \, \ex[\max(X_{T-t},0)] \times f_{M_t}(0).
\label{eq:wh9}
\end{equation}

\subsection{Processes with no down-jumps}

This is a simple alteration of the previous working, pulling the $w$-contour up instead of down. The result is, with $N_t$ denoting the running minimum of $X_t$,
\begin{equation}
\mbox{ [No down-jumps] }\qquad
f_\tau(t)  = 
f_{N_{T-t}}(0)
\times
\frac{1}{t} \, \ex[\max(-X_t,0)] ,
\label{eq:wh10}
\end{equation}
where $f_{N_{T-t}}(0)$ is understood as $\lim_{y\nearrow0} \pr(N_{T-t}>y)/(-y)$.

\subsection{Arithmetic Brownian motion with drift}

This case is obviously a synthesis of the previous two.
As $L(u)=-\half \sigma^2u^2+\I\mu u$, we must have
\[
u^+_s \cdot u^-_s = s\big/\shalf \sigma^2 
\]
so that the second term in (\ref{eq:wh6}) or (\ref{eq:wh7}) can be written similarly to the first, thereby obtaining
\begin{equation}
\mbox{ [Brownian] }\qquad
f_\tau(t)  = 
\frac{1}{T-t} \,\ex[\max(X_{T-t},0)] \times \frac{1}{t} \,\ex[\max(-X_t,0)] \times \frac{2}{\sigma^2},
\label{eq:wh11}
\end{equation}
an alternative form of (\ref{eq:taudens}). Once the factorisation theorem has been derived, this is a much more elegant derivation than that in the Appendix.

%\subsection{$\alpha$-stable}

%The second factor is $\alpha f_{X_t}(0)$.

\section{Examples}

To compare one process with another we can relate the cumulants of the process to observable parameters---drift, volatility, etc---as follows:
\begin{equation}
\hmu = L'(0)/\I, \qquad \hsigma^2 = -L''(0), \qquad \hkappa = \I L'''(0)\big/\hsigma^3.
\end{equation}
For simplicity we refer to these as the cumulants (though the third is technically a normalised cumulant). Without loss of generality we can rescale the process so that $\hsigma=1$.

The two examples that we are about to consider fall into the category of the Carr-Madan-Yor model in which the L\'evy measure is an exponential multiplied by a power of $x$ (and they appear to be the only two tractable ones): more precisely the L\'evy measure is $\propto x^{\gamma}e^{-x/\xi}$, with $\gamma=0,-\frac{3}{2}$. See for example \cite{Schoutens03} for further details.

%As commented above, Unlike Madan \& Schoutens we obtain explicit results rather than leaving the inverse transforms to numerical procedures.

\subsection{Exponential down-jumps plus drift (no diffusion)}

Let $(X_t)$ have down-jumps arising as a Poisson process of rate $\lambda$, with the jumps independently exponentially-distributed of mean $\xi$. A constant drift $\mu>0$ is added, so that $\ex[X_t]=(\mu-\lambda \xi)t$.
Then
\[
L(u) = \I \mu u - \frac{\I  \lambda\xi  u}{1 + \I \xi u}.
\]
Note that $\xi>0$: that the jumps are downwards is taken care of by the signs.
Graphically this process takes the form of a `saw-tooth', with the constant-rate rise regulated by random downward jumps (see the further discussion in Section \ref{sec:conc}).
The cumulants are
\begin{equation}
\hmu = \mu-\lambda\xi, \qquad 
\hsigma^2 = 2\lambda \xi^2, \qquad
\hkappa = -3\big/\sqrt{2\lambda}
\label{eq:edjparam}
\end{equation}
and the  L\'evy-Khinchin representation is
\[
L(u) = \I \mu u  + \int_0^\infty (\lambda/\xi) e^{-x/\xi} (e^{-\I x u}-1) \, dx.
\]
The Wiener-Hopf factors are
\[
\varphi^+_s(u) = \frac{u^+_s}{u^+_s-u}, \qquad \varphi^-_s(u) = \frac{s(1+\I \xi u) }{\mu \xi u^+_s (u^-_s - u)}
\]
where 
\[
u^\pm_s = \frac{\mu - \lambda \xi - \xi s \mp \sqrt{(\mu-\lambda\xi-\xi s)^2+4\mu \xi s} }{2\mu\xi} \, \I.
\]
Viewed as a function of $s$, this is analytic in the complex plane cut along the slit
\[
\mathcal{S} = \big[{-\xi}\inv(\sqrt{\mu}+\sqrt{\lambda\xi})^2, {-\xi}\inv(\sqrt{\mu}-\sqrt{\lambda\xi})^2 \big].
\]
It is marginally easier to deal with the integrals in (\ref{eq:wh6}) than in (\ref{eq:wh7}).
Each is dealt with by collapsing the integral around the branch cut $\mathcal{S}$. Care needs to be taken over the sign of the square root, as the integrand may have a simple pole at $s=0$ depending on the sign of $\mu-\lambda\xi$, and it also affects the behaviour at $\pm\I\infty$.
It is convenient to introduce the function
\[
\mathcal{W}(a,b;t) = \frac{2}{\pi(\sqrt{b}-\sqrt{a})^2} \int_a^b \frac{\sqrt{(b-z)(z-a)}}{z} e^{-zt} \, dz, \qquad 0\le a<b;
\]
note that $\mathcal{W}(a,b;0)=1$ and that its derivative is related to the modified Bessel function $I_1$ by
\[
\deriv{}{t}\, \mathcal{W}(a,b;t) =
-\biggr(\frac{\sqrt{b}+\sqrt{a}}{\sqrt{b}-\sqrt{a}}\biggr) \frac{e^{-(a+b)t/2}}{t} I_1\biggr(\frac{(b-a)t}{2}\biggr).
% -\biggr(\frac{\sqrt{a}+\sqrt{b}}{2}\biggr)^2 e^{-at} \textstyle M\big(\frac{3}{2},3,-(b-a)t\big) .
\]
Then the first integral in (\ref{eq:wh6}), i.e.\ the first factor in the decomposition of the density function, is
\[
\max(\mu-\lambda\xi,0)
+
\frac{1}{2\pi\I} \oint_\mathcal{S} \frac{ \mu - \lambda\xi  - \xi s + \sqrt{(\mu-\lambda\xi-\xi s)^2+4\mu \xi s} }{2 s} \, e^{s(T-t)} \, ds  
\]
which can be written
\begin{equation}
\left.
\begin{array}{rl}
\mu-\lambda\xi + \lambda\xi\, \mathcal{W}\big(\xi\inv(\sqrt{\mu}-\sqrt{\lambda\xi})^2,\xi\inv(\sqrt{\mu}+\sqrt{\lambda\xi})^2;T-t\big) , & \mu\ge \lambda\xi \\[6pt]
 \mu\, \mathcal{W}\big(\xi\inv(\sqrt{\mu}-\sqrt{\lambda\xi})^2,\xi\inv(\sqrt{\mu}+\sqrt{\lambda\xi})^2;T-t\big) , & \mu\le \lambda\xi 
\end{array}
\right\}.
\label{eq:edj1}
\end{equation}
The second integral, i.e.\ the second factor, is
\[
\mu\inv \delta(t) +
\max(\lambda\xi-\mu,0) +
\frac{1}{2\pi\I} \oint_\mathcal{S} \frac{- \mu + \xi \lambda + \xi s + \sqrt{(\mu-\lambda\xi-\xi s)^2+4\mu \xi s} }{2\mu \xi s} \, e^{st} \, ds  .
\]
(The delta-function comes from the fact that the integrand does not decay at $\pm\I\infty$. If one subtracts $\mu\inv$ from the integrand then it is $o(1)$ in that limit and Jordan's lemma can be applied.)
This can be written
\begin{equation}
\left.
\begin{array}{rl}
\displaystyle
\frac{\lambda}{\mu} \, \mathcal{W}\big(\xi\inv(\sqrt{\mu}-\sqrt{\lambda\xi})^2,\xi\inv(\sqrt{\mu}+\sqrt{\lambda\xi})^2;t\big) , & \mu\ge \lambda\xi \\
\displaystyle
\frac{\lambda}{\mu}-\frac{1}{\xi} + \frac{1}{\xi}\, \mathcal{W}\big(\xi\inv(\sqrt{\mu}-\sqrt{\lambda\xi})^2,\xi\inv(\sqrt{\mu}+\sqrt{\lambda\xi})^2;t\big) , & \mu\le \lambda\xi 
\end{array}
\right\}
+ \mu\inv\delta(t) .
\label{eq:edj2}
\end{equation}

\begin{figure}[!htbp]
\begin{center}\begin{tabular}{rl}
(i)&
\scalebox{0.8}{\includegraphics{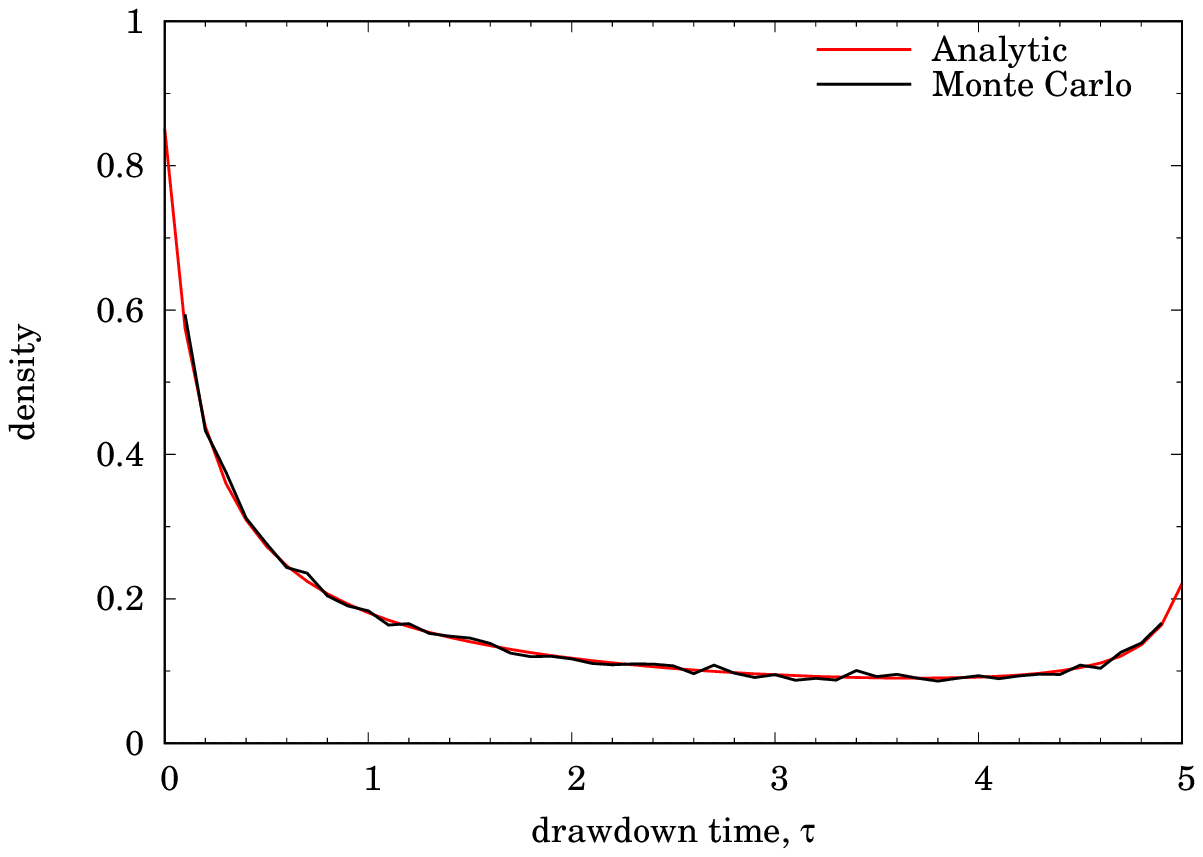}}
\\
(ii)&
\scalebox{0.8}{\includegraphics{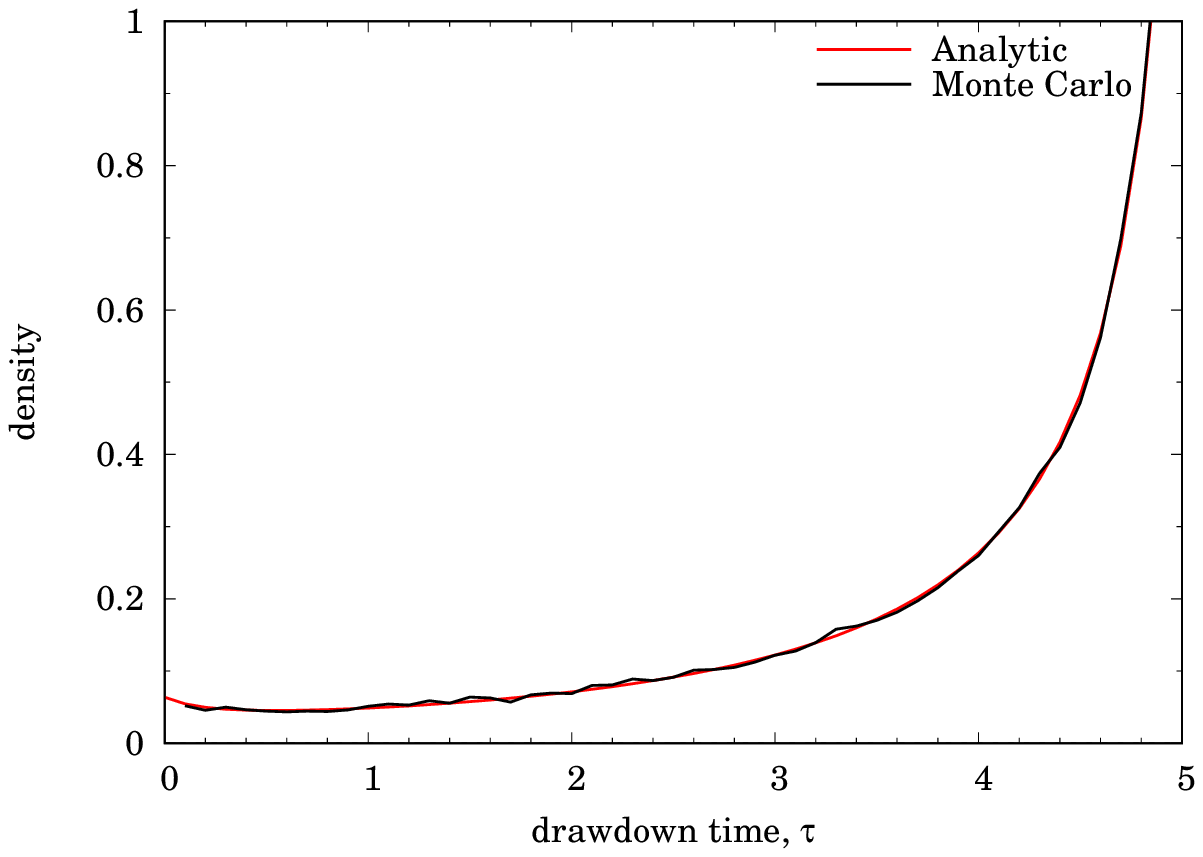}}
\end{tabular}\end{center}
\caption{
Simulation of drawdown time ($10^6$ simulations, bin width 0.1) vs.\ analytical, for process with exponential down-jumps plus drift. Parameters: (i) $\mu=0.6$, $\lambda=4$, $\xi=0.125$; (ii) $\mu=0.3$, $\lambda=4$, $\xi=0.125$; in both cases $T=5$. 
$\pr(\tau=0)=0.213,0.0159$ respectively (not shown in graphs).
}
\label{fig:1}
\end{figure}

\begin{figure}[!htbp]
\centerline{\scalebox{0.8}{\includegraphics{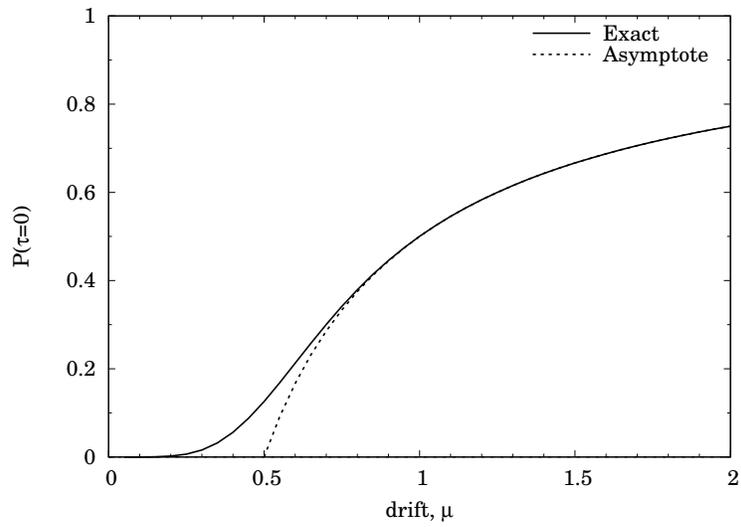}}}
\caption{
Probability of being `at maximum' at time $T$, for the exponential down-jumps plus drift process, as $\mu$ varies: eq.(\ref{eq:pt0},\ref{eq:pt0_approx}) compared.
Parameters: $\lambda=4$, $\xi=0.125$, $T=5$.
}
\label{fig:2}
\end{figure}

Multiplied together, the two expressions (\ref{eq:edj1},\ref{eq:edj2}) give an exact expression for the density function (\ref{eq:wh6}).
Figure~\ref{fig:1} shows a comparison with Monte Carlo.
In addition, it is immediate from (\ref{eq:edj1},\ref{eq:edj2}) that the probability of being `at maximum' at time $T$ is
\begin{equation}
\pr(\tau=0) = \max(1-\lambda\xi/\mu,0) + \min(1,\lambda\xi/\mu)\, \mathcal{W}\big(\xi\inv(\sqrt{\mu}-\sqrt{\lambda\xi})^2,\xi\inv(\sqrt{\mu}+\sqrt{\lambda\xi})^2;T\big) .
\label{eq:pt0}
\end{equation}
This is an increasing function of $\mu$, and as $\mu\searrow0$ it tends to $e^{-\lambda T}$, the probability that no jump occurred, as is intuitively clear.
As $\mu\to\infty$, its behaviour is
\begin{equation}
\pr(\tau=0) \sim 1- \frac{\lambda\xi}{\mu} (1-e^{-\mu T/\xi}).
\label{eq:pt0_approx}
\end{equation}
However $\pr(\tau=T)=0$ because for the maximum to be at time zero the process would have to jump exactly then (an event of measure zero). See Figure~\ref{fig:2} for a plot of this approximation compared with the exact result, in one case.

The Brownian motion limit ($\hkappa\to0$) is reasonably simple to verify, following a careful analysis of $\mathcal{W}$ in that limit.
It is also worth pointing out that when $\mu<0$ none of the results can be used because the Wiener-Hopf factorisation is incorrect at the outset. In fact $\varphi^+_s\equiv1$ (formally, set $u^+_s=-\I C$ and let $C\to\infty$) and so $f_\tau(t)=\delta(T-t)$, obviously, as $(X_t)$ is monotone decreasing when $\mu<0$.

\subsection{Downward Inverse Gaussian process plus drift}

If we take a deterministic positive drift and subtract an Inverse Gaussian process \cite{Schoutens03}, we have
\[
L(u) = \I \mu u +\alpha (1-\sqrt{1+2\beta \I u})
\]
where the square root has positive real part for $u\in \C \setminus [\I/2\beta, \I\infty)$. 
The cumulants are
\[
\hmu = \mu-\alpha\beta, \qquad
\hsigma^2 = \alpha\beta^2, \qquad
\hkappa = -3\big/\sqrt{\alpha}
\]
and the L\'evy-Khinchin representation is
\[
L(u) = \I \mu u + \int_0^{\infty}  \frac{\alpha\sqrt{2\beta} e^{-x/2\beta}}{\sqrt{\pi x^3}} (e^{-\I ux}-1)\, dx.
\]

The Wiener-Hopf factorisation is
\[
\varphi^+_s(u) = \frac{ u^+_s }{u^+_s - u} , \qquad
\varphi^-_s(u) = \frac{ s}{\I \mu u^+_s } \frac{\mu\sqrt{1+2\beta\I u} + \alpha\beta +\sqrt{(\alpha\beta-\mu)^2 + 2\beta\mu s} }{ \mu\sqrt{1+2\beta\I u} - \alpha\beta +\sqrt{(\alpha\beta-\mu)^2 + 2\beta\mu s} }
\]
where the poles are at
\[
u^\pm_s  = \frac{-\alpha^2\beta +\mu\alpha - \mu s  \mp \alpha \sqrt{(\alpha\beta-\mu)^2 + 2\beta\mu s}}{\mu^2} \, \I.
\]
Incidentally in $\varphi^-$ the numerator does not have a zero in the domain of definition, i.e.\ the cut plane: the zero that it does have is on the other branch of the Riemann surface (of the function $u\mapsto\sqrt{1+2\beta\I u}$). Also the denominator does not explicitly invoke $u^-_s$, though of course it could be rewritten.

Performing the Bromwich integrals in (\ref{eq:wh6}) gives the following: the first factor is 
\begin{equation}
(\mu+\alpha \beta) e^{2\alpha (T-t)} \Phi\big({-(\alpha \beta+\mu)} \sqrt{(T-t)/\beta\mu} \big) + (\mu-\alpha\beta) \Phi\big((\mu-\alpha \beta) \sqrt{(T-t)/\beta\mu} \big)
\end{equation}
and the second factor is
\begin{equation}
\frac{1}{\mu} \delta(t) + \frac{2\alpha \sqrt{\beta\mu/t}}{\mu^2} \, \mathcal{C}\big( (\alpha\beta-\mu)\sqrt{t/\beta\mu} \big)
\end{equation}
with $\mathcal{C}(z)$ as defined earlier.
As with the `drift plus exponential down-jumps' model the probability of being at maximum at time $T$ is not zero, and in view of the above it is seen to be
\begin{equation}
(1+\alpha \beta/\mu) e^{2\alpha T} \Phi\big({-(\alpha \beta+\mu)} \sqrt{T/\beta\mu} \big) + (1-\alpha\beta/\mu) \Phi\big((\mu-\alpha \beta) \sqrt{T/\beta\mu} \big).
\end{equation}

\notthis{

\subsection{Downward Reciprocal Inverse Gaussian process with drift}

If we take a deterministic positive drift and subtract an IG process, we have
\[
L(u) = \I \mu u + \alpha \big( (1+2\beta \I u)^{-1/2} -1 \big) 
\]
where the square root has positive real part for $u\in \C \setminus [\I/2\beta, \I\infty)$. 
The cumulants are
\[
\hmu = \mu-\alpha\beta, \qquad
\hsigma^2 = 3 \alpha\beta^2, \qquad
\hkappa = -5\big/6\sqrt{3\alpha}
\]
and the L\'evy-Khinchin representation is ($C$ to be found ***)
\[
L(u) = \I \mu u + \int_0^{\infty}  \frac{C e^{-x/2\beta}}{\sqrt{\pi x}} (e^{-\I ux}-1)\, dx.
\]

Problem here is that cannot get the position of the pole in closed form: needs a cubic equation to be solved.
}

\subsection{Graphical comparison and discussion}

As suggested earlier, we wish to study the effect of generalising  the arithmetic Brownian motion, and the natural way to do this is to use the $(\hmu,\hsigma,\hkappa)$ parameterisation introduced earlier. By rescaling we can choose to fix $T$ and $\hsigma$, so that processes are classified by $\hmu$ and $\hkappa$.

Rather than plotting the density, which does not allow the delta-function component to be represented in a meaningful way, we instead plot the cumulative as obtained by numerically integrating the density using the trapezium rule (except for the arithmetic Brownian motion case, $\hkappa=0$, where we can use (\ref{eq:taucumul})). Without loss of generality we can consider $\hkappa<0$, as the opposite case can be considered by point symmetry ($t\mapsto T-t$, $X\mapsto-X$).

\begin{figure}[!htbp]
\begin{center}\begin{tabular}{rl}
(i)&
\scalebox{0.8}{\includegraphics{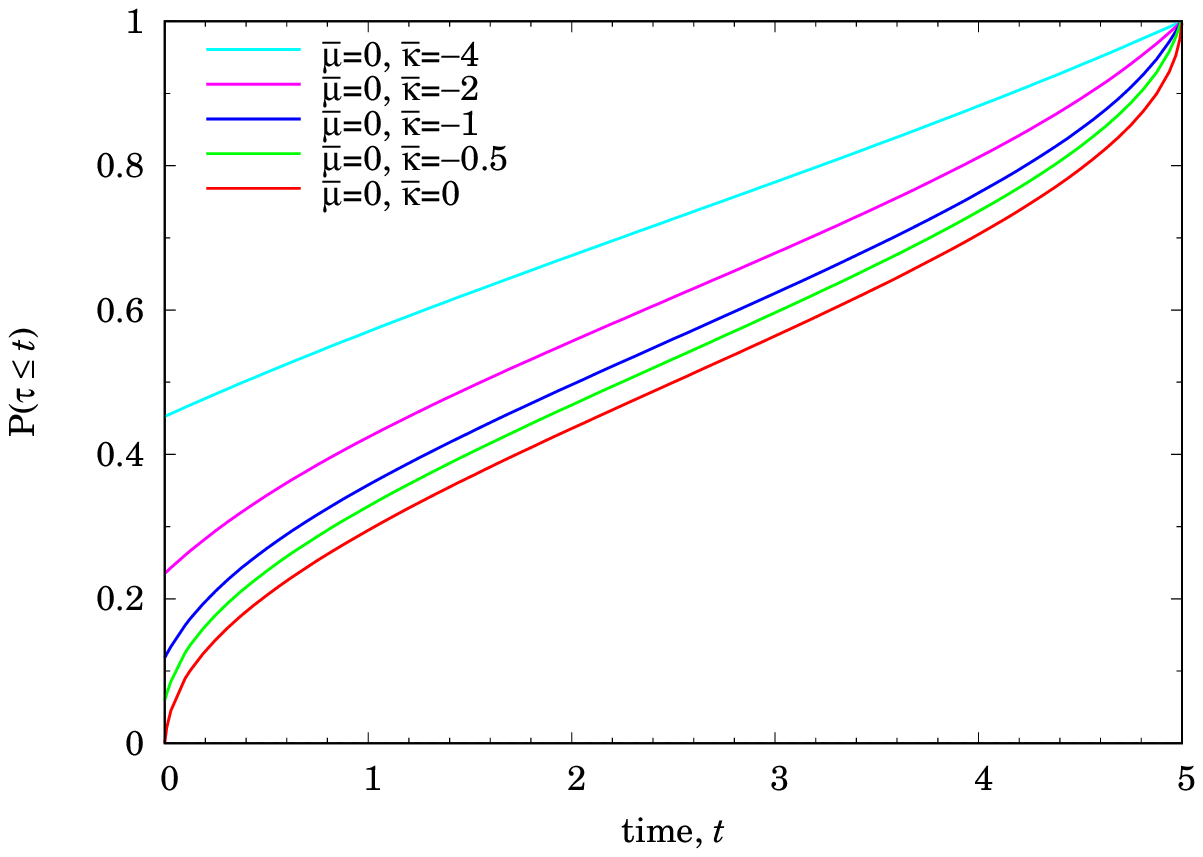}} \\
(ii)&
\scalebox{0.8}{\includegraphics{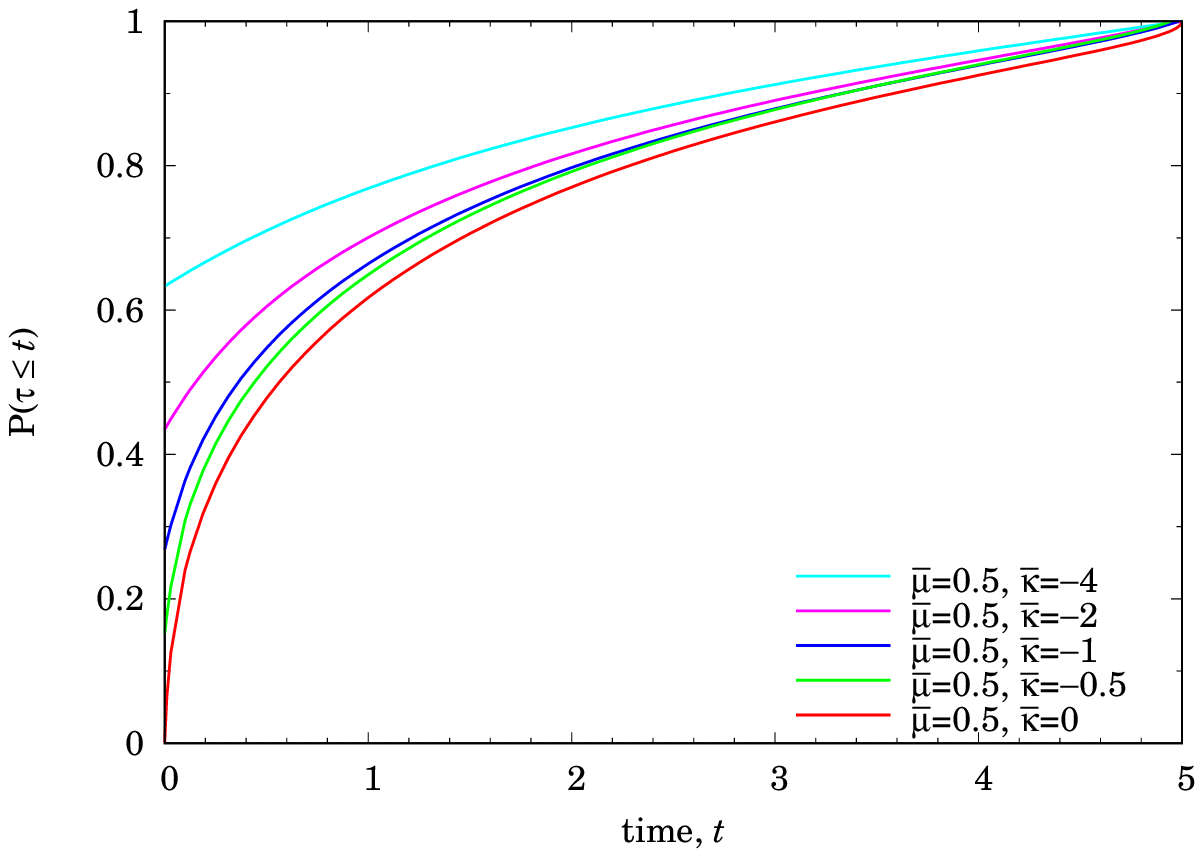}} \\
(iii)&
\scalebox{0.8}{\includegraphics{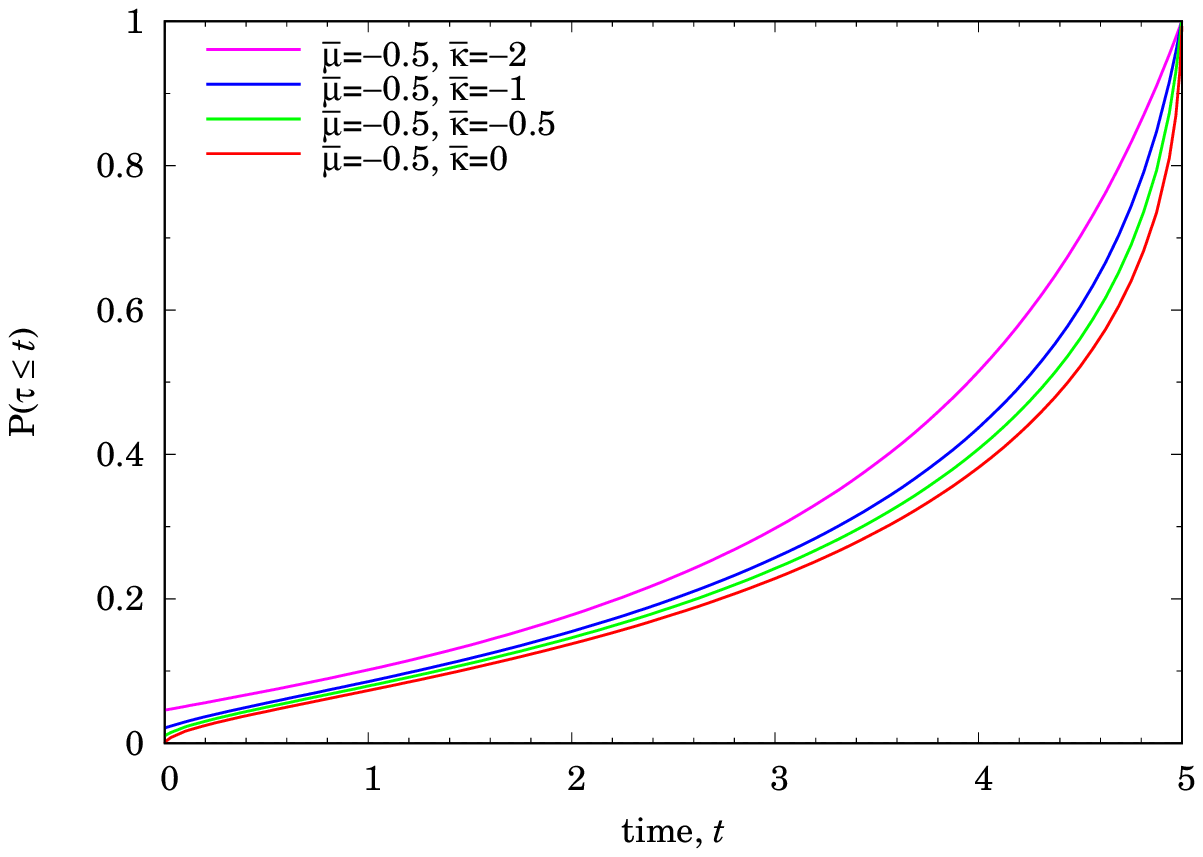}}
\end{tabular}\end{center}
\caption{
Distribution of time since maximum for exponential down-jumps process with drift: (i) $\hmu=0$, (ii) $\hmu=+0.5$, (iii) $\hmu=-0.5$. In each $\hkappa$ is varied, and $\hsigma$ is fixed at 1; $T=5$ throughout.
}
\label{fig:edj}
\end{figure}

\begin{figure}[!htbp]
\begin{center}\begin{tabular}{rl}
(i)&
\scalebox{0.8}{\includegraphics{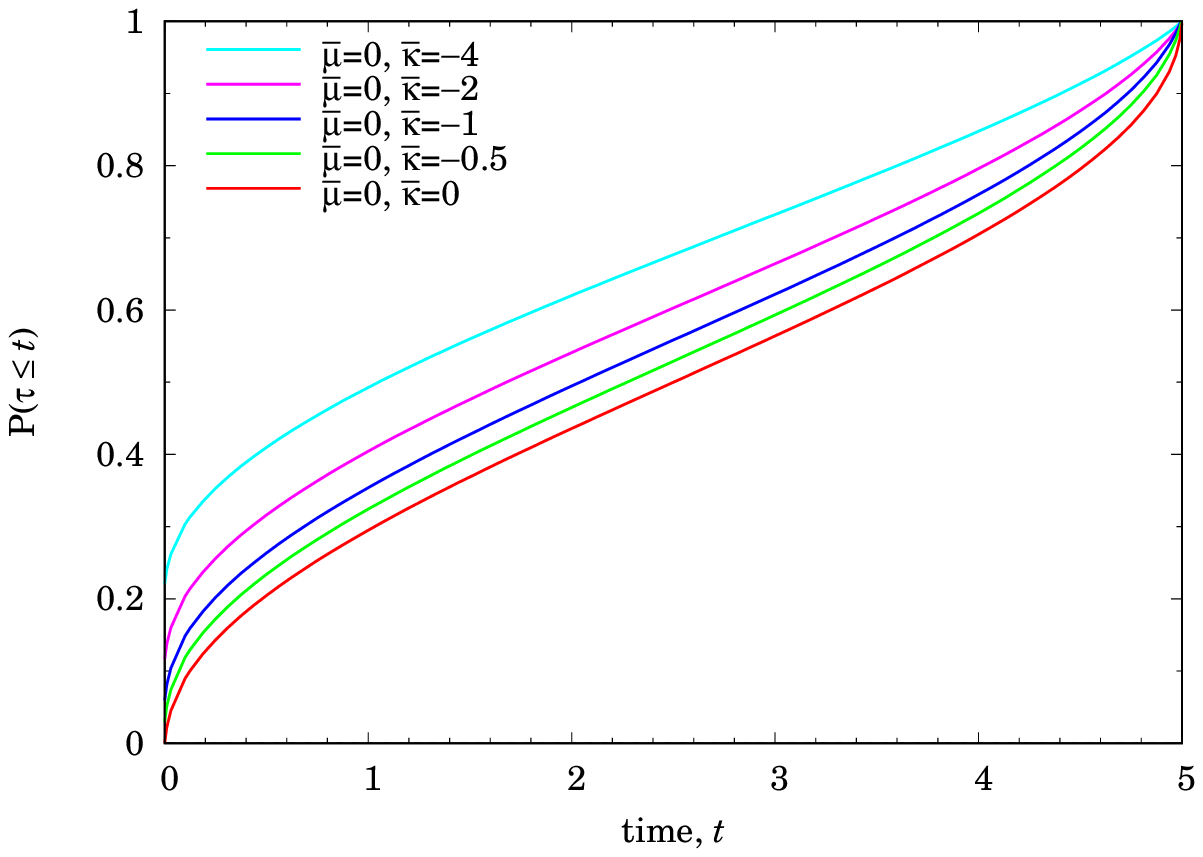}} \\
(ii)&
\scalebox{0.8}{\includegraphics{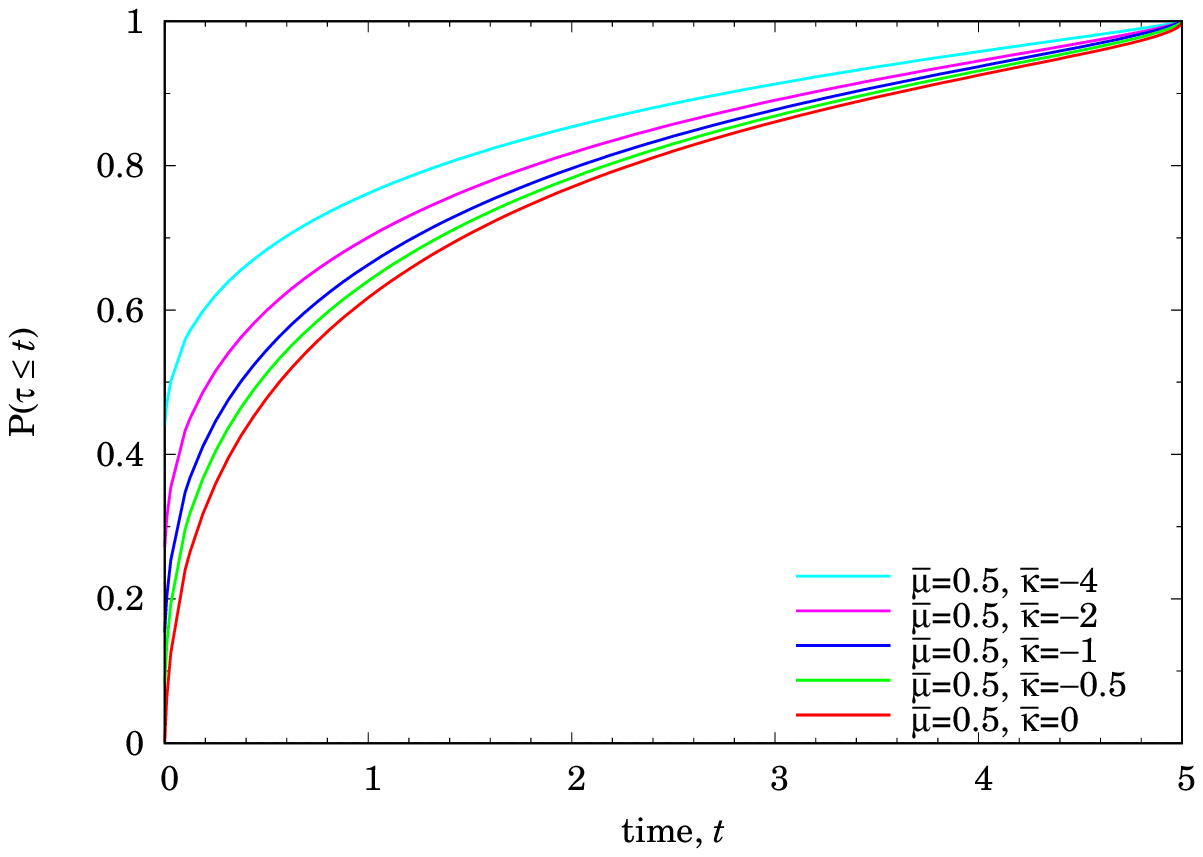}} \\
(iii)&
\scalebox{0.8}{\includegraphics{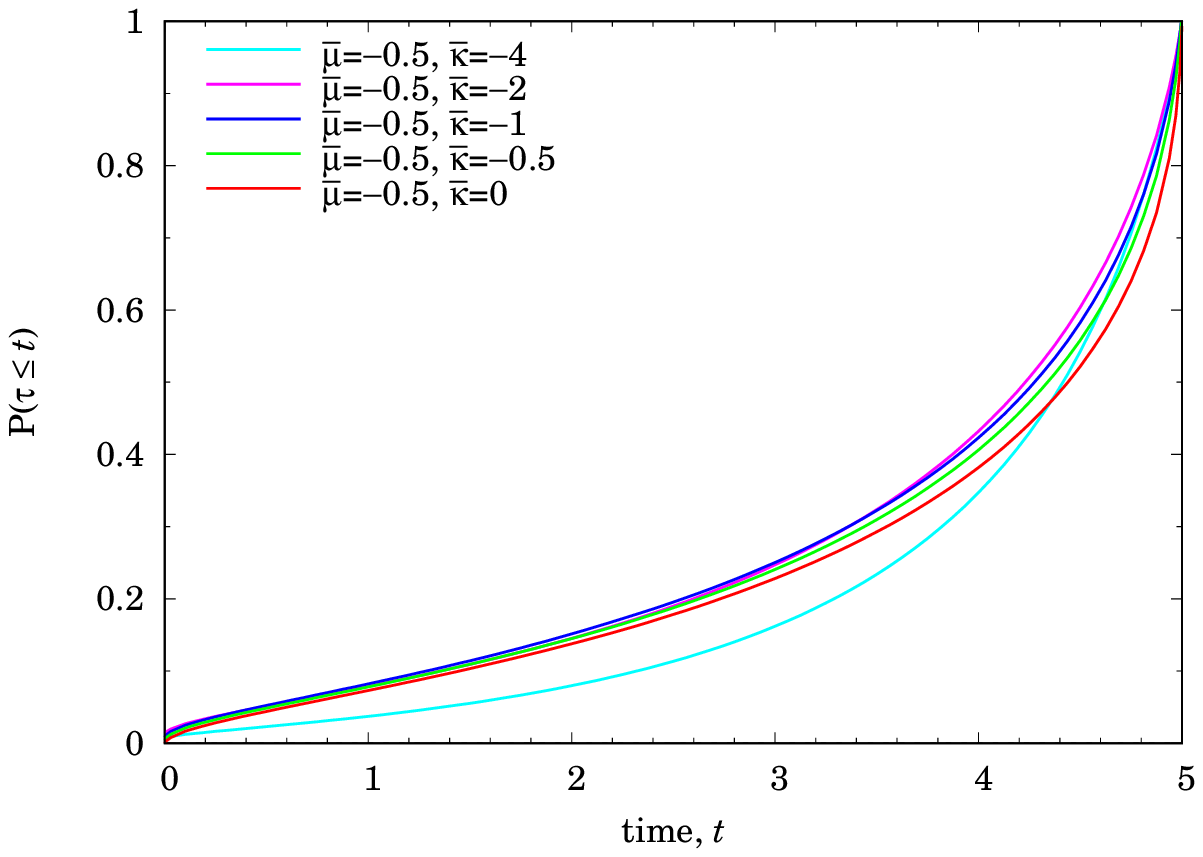}}
\end{tabular}\end{center}
\caption{
Distribution of time since maximum for downward Inverse Gaussian process with drift: (i) $\hmu=0$, (ii) $\hmu=+0.5$, (iii) $\hmu=-0.5$. In each $\hkappa$ is varied, and $\hsigma$ is fixed at 1; $T=5$ throughout.
}
\label{fig:ig}
\end{figure}

Figure~\ref{fig:edj} shows the effect of varying $\hkappa$ for different values of $\hmu$ for the exponential down-jumps model. Figure~\ref{fig:ig} repeats this for the Inverse Gaussian.

When $\hkappa=-4$, with the exponential down-jumps model, the drift $\mu$ becomes negative if $\hmu=-0.5$, and in that case the process is monotone decreasing, and $\tau=T$ with probability 1. We make the point that the graph is not a continuous function of $\mu$ as $\mu$ crosses zero.

An obvious departure from the arithmetic Brownian motion in these models is that when $\hkappa<0$ there is positive probability that $\tau=0$. This would not be observed if, for example, we used a diffusion \emph{and} down-jumps as well. However, processes of that sort do not appear to be as analytically tractable as the two we have shown here. In the case of diffusion plus exponential down-jumps, the L\'evy generator is
\[
L(u) = \I \mu u - \shalf \sigma^2 u^2 - \frac{\I  \lambda\xi  u}{1 + \I \xi u}
\]
and now $u^+_s$ can only be obtained by factorising a cubic. The resulting inverse Laplace transforms in (\ref{eq:wh6}) then have to be done numerically.

It is apparent that when the total drift $\hmu$ is negative, a negative skewness $\hkappa$ has less of an effect than when  $\hmu$ is positive.

\section{Conclusions and final remarks}
\label{sec:conc}

This paper has considered the time since maximum over some fixed observation interval for a L\'evy process, as a generalisation of Brownian motion. It relies on what we are calling the factorisation theorem, which states that for a completely asymmetric L\'evy process, i.e.\ one in which the jumps are not in both directions, the density of the time since maximum splits as a product of two parts, one dependent on $t$ and the other on $T-t$, as expressed variously in eq.~(\ref{eq:wh6},\ref{eq:wh7},\ref{eq:wh9},\ref{eq:wh10},\ref{eq:wh11}).
The computation of these factors depends on the Wiener-Hopf factorisation, which can be computed for this type of process as one of the factors is a simple pole. For reasons of symmetry, the same statements are true of the time to maximum, i.e.\ measured from the origin.

One of the examples that we have solved explicitly is that of exponential down jumps in the presence of a constant positive drift. This three parameter model exhibits rich behaviour and has Brownian motion as one of its limiting forms. In that sense it provides a model case to study deviations from conventional Brownian motion. Moreover, the `saw-tooth' profile is of interest in its own right, with potential application in queuing systems including traffic flow and inventory processes \cite{Asmussen03}. In such cases, the time to reach maximum over a given observation window is often a key variable of interest.

It is worth mentioning a connection between our work and that of \cite{Challet17} in relation to investment management. Our conclusions suggest that for processes with the same Sharpe ratio $\hmu/\hsigma$, the effect of making the skewness negative is to increase the probability of short or even zero-length drawdowns. This is unsurprising as the process tends to crawl upwards and occasionally jump downwards. Despite the presence of a greater number of `upper records', in the  terminology of \cite[\S1]{Challet17}, which is superficially comforting, such an investment is clearly \emph{worse} than one of the same Sharpe ratio but zero skew. 
Asset classes in which this may occur are illiquid credit strategies and insurance: in both, one receives a steady coupon but runs the risk of occasional large drawdowns.

As a general point, if $\tau$ is some random variable that lies in the interval $[0,T]$ a.s., then the factorisation theorem is saying that
\[
f_\tau(t) = \biggr(\frac{1}{2\pi\I} \invintip \frac{1}{G(s)} e^{(T-t)s} \, ds  \biggr)
\biggr(\frac{1}{2\pi\I} \invintip \frac{G(s)}{s} e^{ts} \, ds \biggr)
\]
satisfies $\int_0^T f_\tau(t)\, dt=1 $ for all $T$, as is immediate from the convolution theorem, and it is in effect what is going on here. Therefore it is, at least in principle, a way of factorising the density of any such random variable. Some restrictions on $G$ are needed to ensure that the RHS be a valid p.d.f.: at least, it needs to be analytic and free from zeros in the right half-plane. Whether this construction is more generally applicable to other problems in stochastic processes is a matter for further research.
One possible area would be to look at constrained processes, extending what is known for the Brownian motion case \cite{Majumdar08a}.

When the process does have both up- and down-jumps, the factorisation theorem no longer holds, on account of there being multiple singularities in the upper and lower half planes. In that case the result has to be replaced by a sum of products, if the singularities are poles, or an integral (continuous sum)  in the general case.

\notthis{

\section{Computation}

Computational note on $\mathcal{W}$.
When $(b-a)t$ is large, the majority of the variation in the integrand comes from the exponential term. We therefore write
\[
\mathcal{W}(a,b;t) =
\frac{2e^{-at}}{\pi^{1/2}(\sqrt{b}-\sqrt{a})^2t} \int_0^\infty \left\{ \frac{u\sqrt{(b-a)t-u}}{u+at} \, \indic{u<(b-a)t}\right\} \frac{e^{-u}}{\sqrt{\pi u}} \, du
\]
and use the Gauss-Laguerre quadrature formula to evaluate the integral. (Extracting a factor of $u^{-1/2}$ causes the bracketed term to vary more slowly at the origin, particularly for $at\ll1$.)
When $(b-a)t$ is small, we write instead
\[
\mathcal{W}(a,b;t) = (\sqrt{b}+\sqrt{a})^2 e^{-at} \int_0^1  \left\{ \frac{ue^{-(b-a)tu}}{(b-a)u+a} \right\} \frac{u^{-1/2}(1-u)^{1/2}}{\pi/2} \, du
\]
which is then done by Gauss-Jacobi. A practical proposition is to use the first method for $(b-a)t\gtrsim5$ and the second for $(b-a)t\lesssim 5$.

} % end notthis

\appendix

\section{Appendix: Brownian motion derivation and formulary}

We start with some useful results that pertain to integrals of the functions $\phi$, $\Phi$ and move on to some identities concerning the bivariate Normal integral. After that we derive (\ref{eq:taudens},\ref{eq:taucumul}).

\subsection{Integrals of products of Normal distribution functions}

Let $a_1$, $b_1$, $a_2$, $b_2$ be real constants, and write $b^* = \sqrt{b_1^2+b_2^2}$, $a^*=(a_1b_1+a_2b_2)/b^*$.
Then:
\begin{eqnarray*}
\int_0^\infty \phi(a_1-b_1x) \phi(a_2-b_2x)\, dx &=& \frac{1}{b^*} \phi\!\left( \frac{a_1b_2-a_2b_1}{b^*} \right) \Phi(a^*), \qquad b_1>0, \quad b_2>0 \\
\int_0^\infty \phi(a_1-b_1x) \Phi(a_2-b_2x)\, dx &=& \frac{1}{b_1} \Phi_2\!\left( a_1, \frac{a_2b_1-a_1b_2}{b^*} ; \frac{-b_2}{b^*}\right), \qquad b_1>0 \\
\int_0^\infty \phi(a_1+b_1x) \Phi(a_2-b_2x)\, dx &=& \frac{1}{b_1} \Phi_2\!\left(-a_1, \frac{a_2b_1+a_1b_2}{b^*} ; \frac{-b_2}{b^*}\right), \qquad b_1>0
\end{eqnarray*}
The first is immediate.
The second and third follow from the identity
\[
\Phi_2(x,y;\rho) = \int_{-\infty}^x \Phi\!\biggr(\frac{y-\rho\xi}{\sqrt{1-\rho^2}}\biggr) \phi(\xi)\, d\xi
\]
which is derived thus: if $(X,Y)$ is unit bivariate Normal with correlation $\rho$ then we can write $Y=\rho X+\sqrt{1-\rho^2}Z$ with $Z,X$ independent and $Z$ standard Normal, condition on $X$ and then integrate out:
\[
\pr(X<x \, \cap\, Y<y) = \ex\left[ \pr\!\left(Z<\frac{y-\rho X}{1-\rho^2} \,\Big|\, X \right) \indic{X<x} \right].
\]
Finally we note the quadrant law \cite{Abramowitz64}:
\[
\Phi_2(-x,-y;\rho) = \Phi_2(x,y;\rho) - \Phi(x)-\Phi(y)+1
\]

\subsection{Reciprocity law for bivariate Normal integral}

If $\rho\ge 0$ and $\rho^*=\sqrt{1-\rho^2}$, then
\begin{equation}
\left.
\begin{array}{rcl}
\Phi_2(-x,\rho x;-\rho) + \Phi_2(x,-\rho^*x;-\rho^*) &=& \Phi(\rho x) \Phi(-\rho^*x) \\[6pt]
\Phi_2(-x,-\rho x;\rho) + \Phi_2(x,\rho^*x;\rho^*) &=&  1 - \Phi(\rho x) \Phi(-\rho^*x)
\end{array}
\right\}
\label{eq:recip}
\end{equation}
Proof.
In the LHS of either expression, write $\Phi_2$ as a double integral over the bivariate Normal density and then differentiate with respect to $x$. A little algebra allows this resulting expression to be manipulated into $\rho\phi(\rho x)\Phi(-\rho^*x)-\rho^*\Phi(\rho x)\phi(\rho^*x)$ in the first case, and $(-1\times\mbox{this})$ in the second. When $x=0$ the identity is clear from the orthant law
\begin{equation}
\Phi_2(0,0;\rho) = \frac{1}{4} + \frac{1}{2\pi}\arcsin \rho
\label{eq:orthant}
\end{equation}
and as $\arcsin \rho + \arcsin \rho^* = \frac{\pi}{2}$, 
this completes the proof. $\Box$

A use for this law is that the tetrachoric series expansion of $\Phi_2$ converges very slowly as $\rho\to1$, but the reciprocity law shows that one only need evaluate $\Phi_2(x,\rho x;\rho)$ or $\Phi_2(-x,\rho x;-\rho)$ for $\rho\le \frac{1}{\sqrt{2}}$.
From this we deduce the previously-mentioned
\begin{equation}
\Phi_2\! \left( x,\frac{x}{\sqrt{2}} ; \frac{1}{\sqrt{2}}  \right) = \frac{1}{2} \left( \Phi\!\left( \frac{x}{\sqrt{2}}\right)^2 + \Phi(x) \right) .
\label{eq:recip2}
\end{equation}

\subsection{Derivation of (\ref{eq:taudens}) and (\ref{eq:taucumul})}

The joint density of $(X_t,M_t)$ is (e.g.~\cite{Borodin02})
\begin{equation}
\frac{2(2y-x)}{\sqrt{2\pi \sigma^6 t^3}} \exp\left(  \frac{-(2y-x+\mu t)^2 + 4\mu yt}{2\sigma^2 t} \right) \indic{x<y} \indic{0<y} 
\end{equation}
while the marginal density of $M_t$ is
\begin{equation}
f_{M_t}(y) = \frac{2}{\sigma\sqrt{t}} \phi \! \left( \frac{y-\mu t}{\sigma\sqrt{t}} \right) - \frac{2\mu}{\sigma^2} e^{2\mu y/\sigma^2} \Phi \! \left( \frac{-y-\mu t}{\sigma\sqrt{t}} \right)\indic{0<y} .
\end{equation}
and its cumulative is
\begin{equation}
\pr(M_t>y) = \Phi \! \left( \frac{\mu t - y}{\sigma\sqrt{t}} \right) +  e^{2\mu y/\sigma^2} \Phi \! \left( \frac{-\mu t -y}{\sigma\sqrt{t}} \right).
\label{eq:cumulM}
\end{equation}

Viewed from time $T$, the probability of that the maximum occurred at most a time $T-t$ ago, i.e.\ on or after time $t$, is just $\pr(M_T>M_t)$.
To find this, we first condition on $(X_t,M_t)$ and then integrate out. As the process is Markovian,
$\pr(M_T>y \cdl X_t=x)$ is obtained from (\ref{eq:cumulM}) simply by replacing $t$ with $T-t$ and $y$ with $y-x$. We then integrate $x$ and $y$ out. Accordingly
\begin{eqnarray*}
\lefteqn{\pr(\tau<T-t) = \pr(M_T>M_t) = \int_0^\infty \int_{-\infty}^y
\frac{2(2y-x)}{\sqrt{2\pi \sigma^6 t^3}} 
\exp\left(  \frac{-(2y-x+\mu t)^2 + 4\mu yt}{2\sigma^2 t} \right)
} \\
&& \times 
\left[ \Phi \! \left( \frac{\mu (T-t) - (y-x)}{\sigma\sqrt{T-t}} \right) +  e^{2\mu (y-x)/\sigma^2} \Phi \! \left( \frac{-\mu (T-t) -(y-x)}{\sigma\sqrt{T-t}} \right) \right]
\, dx \, dy .
\end{eqnarray*}
Thus on replacing $t$ by $T-t$ and changing variable from $x$ to $y-x$ in the inner integral we find
\begin{eqnarray*}
\lefteqn{\hspace{-30mm}\pr(\tau<t) = \int_0^\infty \int_0^\infty
\frac{2(y+x)}{\sqrt{2\pi \sigma^6 (T-t)^3}} 
\exp\left(  \frac{-\big(y+x+\mu (T-t)\big)^2 + 4\mu y(T-t)}{2\sigma^2 (T-t)} \right)
} \\
&& \times 
\left[ \Phi \! \left( \frac{\mu t - x}{\sigma\sqrt{t}} \right) +  e^{2\mu x/\sigma^2} \Phi \! \left( \frac{-\mu t -x}{\sigma\sqrt{t}} \right) \right]
\, dx \, dy 
\end{eqnarray*}
whereupon the $y$-integral can be done to give
\begin{eqnarray}
\lefteqn{\hspace{-30mm}\pr(\tau<t) = \int_0^\infty \left[ \frac{2}{\sigma\sqrt{T-t}} \phi\!\left(\frac{x+\mu(T-t)}{\sigma\sqrt{T-t}}\right) + \frac{2\mu}{\sigma^2} e^{-2\mu x/\sigma^2} \Phi\!\left(\frac{\mu (T-t)-x}{\sigma\sqrt{T-t}}\right) \right]  } \nonumber \\
&& \times 
\left[ \Phi \! \left( \frac{\mu t - x}{\sigma\sqrt{t}} \right) +  e^{2\mu x/\sigma^2} \Phi \! \left( \frac{-\mu t -x}{\sigma\sqrt{t}} \right) \right]
\, dx  
\label{eq:messy1a}
\end{eqnarray}
after which the first square bracket can be integrated to allow the following representation:
\begin{eqnarray}
\lefteqn{\hspace{-30mm}\pr(\tau<t) = - \int_0^\infty \pderiv{}{x} 
\left[ \Phi\!\left(\frac{- \mu (T-t)-x}{\sigma\sqrt{T-t}}\right)  + e^{-2\mu x/\sigma^2}  \Phi\!\left(\frac{\mu(T-t)-x}{\sigma\sqrt{T-t}}\right) \right]  } \nonumber \\
&& \times 
\left[ \Phi \! \left( \frac{\mu t - x}{\sigma\sqrt{t}} \right) +  e^{2\mu x/\sigma^2} \Phi \! \left( \frac{-\mu t -x}{\sigma\sqrt{t}} \right) \right] \label{eq:messy1b}
\, dx  .
\end{eqnarray}
Write $Q(t;\mu)$ for $\pr(\tau<t)$ and integrate (\ref{eq:messy1b}) by parts. Then, as expected,
\[
Q(t;\mu) + Q(T-t;-\mu) = 1
\]
so that density is invariant under the point-symmetry transformation.

The next step is to obtain the density by differentiating (\ref{eq:messy1b}) w.r.t.\ $t$ and then doing the integrals. Eventually (\ref{eq:taudens}) is obtained, but the derivation is laborious; some useful short-cuts are given in \S{A.1}. Interestingly, the bivariate Normal integral appears, but cancels itself out through what we are calling the `reciprocity law' (\ref{eq:recip}), so the end result contains no $\Phi_2$ term.
Unsurprisingly $f_\tau(t)\cdot T$ depends on $t$ through two dimensionless quantities: the ratio $\mu\sqrt{T}/\sigma$ (that is, the mean at time $T$ divided by the standard deviation), and the ratio $t/T$.

We now return to the problem of finding $\pr(\tau<t)$. Return to (\ref{eq:messy1a}) and write it in short as $\int[A+B][C+D]$. Then $AC+AD$ gives
\[
2\Phi_2\!\left(
\frac{-\mu\sqrt{T-t}}{\sigma},\frac{\mu\sqrt{T}}{\sigma};
-\sqrt{\frac{T-t}{T}}
\right) +
2\Phi_2\!\left(
\frac{\mu\sqrt{T-t}}{\sigma},\frac{-\mu\sqrt{T}}{\sigma};
-\sqrt{\frac{T-t}{T}}
\right) 
\]
while $BC$ gives (on integration by parts)
\[
\Phi\!\left(\frac{\mu\sqrt{T-t}}{\sigma} \right)
\Phi\!\left(\frac{\mu\sqrt{t}}{\sigma} \right)
- \Phi_2 \! \left( \frac{-\mu\sqrt{T-t}}{\sigma}, \frac{\mu\sqrt{T}}{\sigma} ; -\sqrt{\frac{T-t}{T}} \right)
- \Phi_2 \! \left( \frac{-\mu\sqrt{t}}{\sigma}, \frac{\mu\sqrt{T}}{\sigma} ; -\sqrt{\frac{t}{T}} \right)
\]
and $BD$ gives (multiply the integrand by unity and integrate by parts)
\[
\frac{2\mu}{\sigma^2} \left\{
 \sigma\sqrt{T-t} \, \phi\!\left(\frac{\mu\sqrt{T-t}}{\sigma}\right) \Phi\!\left(-\frac{\mu\sqrt{t}}{\sigma}\right) 
+
\sigma\sqrt{t} \, \phi\!\left(\frac{\mu\sqrt{t}}{\sigma}\right)  \Phi\!\left(\frac{\mu\sqrt{T-t}}{\sigma}\right) 
- \frac{\sigma\sqrt{T}}{2} \, \phi \left( \frac{\mu\sqrt{T}}{\sigma} \right)
\right.
\]
\[
\left.
\mbox{} + \mu (T-t) \Phi_2\!\left(
\frac{\mu\sqrt{T-t}}{\sigma},\frac{-\mu\sqrt{T}}{\sigma};
-\sqrt{\frac{T-t}{T}}
\right) 
- \mu t \Phi_2\!\left(
\frac{-\mu\sqrt{t}}{\sigma},\frac{\mu\sqrt{T}}{\sigma};
-\sqrt{\frac{t}{T}}
\right) 
\right\}.
\]
The sum of these therefore gives the required cumulative distribution function:
\[
\Phi\!\left(\frac{\mu\sqrt{T}}{\sigma}\right) 
+ 2\Phi_2\!\left(\frac{\mu\sqrt{T-t}}{\sigma},\frac{-\mu\sqrt{T}}{\sigma};-\sqrt{\frac{T-t}{T}}\right) 
- 2\Phi_2\!\left(\frac{-\mu\sqrt{t}}{\sigma},\frac{\mu\sqrt{T}}{\sigma};-\sqrt{\frac{t}{T}}\right) 
\]
\[
\mbox{}+ \frac{2\mu}{\sigma^2} \left\{
 \sigma\sqrt{T-t} \, \phi\!\left(\frac{\mu\sqrt{T-t}}{\sigma}\right) \Phi\!\left(\frac{-\mu\sqrt{t}}{\sigma}\right) 
+
\sigma\sqrt{t} \, \phi\!\left(\frac{\mu\sqrt{t}}{\sigma}\right)  \Phi\!\left(\frac{\mu\sqrt{T-t}}{\sigma}\right) 
- \frac{\sigma\sqrt{T}}{2} \, \phi \left( \frac{\mu\sqrt{T}}{\sigma} \right)
\right.
\]
\[
\left.
\mbox{}+ \mu (T-t) \Phi_2\!\left(
\frac{\mu\sqrt{T-t}}{\sigma},\frac{-\mu\sqrt{T}}{\sigma};
-\sqrt{\frac{T-t}{T}}
\right) 
- \mu t \Phi_2\!\left(
\frac{-\mu\sqrt{t}}{\sigma},\frac{\mu\sqrt{T}}{\sigma};
-\sqrt{\frac{t}{T}}
\right) 
\right\}.
\]
This can be further tidied-up to give (\ref{eq:taucumul}).

\section*{Acknowledgements}

We thank Alex Mijatovi\'c for helpful discussions.

\bibliographystyle{plain}
\bibliography{../phd}

\end{document}